\documentclass[12pt,thmsa]{article}
\usepackage{sw20bams}

\input tcilatex
\begin{document}

\author{Steven Finch}
\title{Perimeter Variance of Uniform Random Triangles }
\date{July 1, 2010}
\maketitle

\begin{abstract}
Let $T$ be a random triangle in a disk $D$ of radius $R$ (meaning that
vertices are independent and uniform in $D$). We determine the bivariate
density for two arbitrary sides $a,b$ of $T$. In particular, we compute that 
$\limfunc{E}(a\,b)=(0.837...)R^2$, which implies that $\limfunc{Var}($%
perimeter$)=(0.649...)R^2$. No closed-form expression for either coefficient
is known. The Catalan numbers also arise here.
\end{abstract}

\footnotetext{
Copyright \copyright\ 2010 by Steven R. Finch. All rights reserved.}Let $%
A,B,C$ denote three independent uniformly distributed points in the disk $%
D=\{(\xi ,\eta ):\xi ^2+\eta ^2\leq R^2\}$. Let $T$ denote the triangle with
sides $a,b,c$ opposite the vertices $A,B,C$. We are interested in the
perimeter $a+b+c$ of triangle $T$. The univariate density $f(x)$ for side $a$
is \cite{Delt, Ham, Lord, Alag, Solo, Dunb, TuFi} 
\[
\begin{array}{ccc}
\dfrac{4x}{\pi R^2}\arccos \left( \dfrac x{2R}\right) -\dfrac{x^2}{\pi R^4}%
\sqrt{4R^2-x^2}, &  & 0<x<2R
\end{array}
\]
and 
\[
\begin{array}{ccc}
\limfunc{E}(a)=\dfrac{128}{45\pi }R=(0.9054147873672267990407609...)R, &  & 
\limfunc{E}(a^2)=R^2.
\end{array}
\]
Clearly 
\[
\limfunc{E}(\text{perimeter})=3\limfunc{E}(a)=\dfrac{128}{15\pi }%
R=(2.7162443621016803971222828...)R 
\]
but to compute $\limfunc{Var}($perimeter$)=\limfunc{E}($perimeter$^2)-%
\limfunc{E}($perimeter$)^2$, we will further need to consider
cross-correlation $\rho $ between sides.

The bivariate density $f(x,y)$ for sides $a,b$ is 
\[
f(x,y)=\left\{ 
\begin{array}{lll}
\varphi (x,y) &  & \text{if }x+y\leq 2R, \\ 
\psi (x,y) &  & \text{if }x+y>2R\text{ and }x\leq 2R
\end{array}
\right. 
\]
when $0\leq y\leq x$ (use symmetry otherwise) where 
\begin{eqnarray*}
\varphi (x,y) &=&\frac{2xy}{\pi R^6}\left\{ -\sqrt{(2R-x-y)(x-y)(2R+x-y)(x+y)%
}\,+\right. \\
&&\ \left. 2(R-y)^2\arccos \left( \frac{x^2-2Ry+y^2}{2x(R-y)}\right)
+2R^2\arccos \left( \frac{x^2+2Ry-y^2}{2Rx}\right) \right\} + \\
&&\ \frac{8xy}{\pi ^2R^6}\dint\limits_{R-y}^Rt\arccos \left( \frac{%
t^2+x^2-R^2}{2tx}\right) \arccos \left( \frac{t^2+y^2-R^2}{2ty}\right) dt,
\end{eqnarray*}
\[
\psi (x,y)=\frac{8xy}{\pi ^2R^6}\dint\limits_{x-R}^Rt\arccos \left( \frac{%
t^2+x^2-R^2}{2tx}\right) \arccos \left( \frac{t^2+y^2-R^2}{2ty}\right) dt. 
\]
It follows by numerical integration that 
\[
\limfunc{E}(a\,b)=(0.8378520652962219016710654...)R^2 
\]
hence 
\[
\rho (a,b)=\dfrac{\limfunc{E}(a\,b)-\limfunc{E}(a)\limfunc{E}(b)}{\sqrt{%
\limfunc{Var}(a)\limfunc{Var}(b)}}=0.1002980835659001715822627..., 
\]
\[
\limfunc{E}(\text{perimeter}^2)=3\limfunc{E}(a^2)+6\limfunc{E}%
(a\,b)=(8.0271123917773314100263929...)R^2, 
\]
\[
\limfunc{Var}(\text{perimeter})=(0.6491289571281667551974101...)R^2. 
\]
Exact evaluation of $\limfunc{E}(a\,b)$ remains an open problem. We review
derivation of the univariate case in section 1, imitating the analysis in 
\cite{Parry, PaFi} very closely. (Parry's thesis \cite{Parry} is concerned
with triangles in three-dimensional space; it is surprising that our
two-dimensional analog has not yet been examined.) The bivariate case is
covered in section 2. An experimental consequence of our work is the formula 
\[
\limfunc{E}(a^2\,b^2)=\frac{13}{12}R^4 
\]
which we prove via a different approach in section 3. Finally, in section 4,
the Catalan numbers from combinatorics appear rather unexpectedly.

\section{Univariate Case}

We omit geometric details, referring to \cite{Parry, PaFi} instead. The
distance $t$ between point $C$ and the origin has density $2t/R^2$ for $%
0<t<R $. Let $f(x\,|\,t)$ be the conditional density for distance $x$
between points $C$ and $B$, given $t$. We will compute the sought-after
density $f(x) $ for side $a$ via 
\[
f(x)=\dint\limits_0^Rf(x\,|\,t)\frac{2t}{R^2}dt. 
\]
There are two subcases.

\subsection{$0<x<R$}

\[
f(x)=\dint\limits_0^{R-x}\frac{2x}{R^2}\frac{2t}{R^2}dt+\dint\limits_{R-x}^R%
\frac{2x}{\pi R^2}\arccos \left( \frac{t^2+x^2-R^2}{2tx}\right) \frac{2t}{R^2%
}dt 
\]
which corresponds to formula (1.12) in Parry's thesis \cite{Parry}. The
arccos term arises since, if the portion of a circle of radius $x$, center $%
C $ contained within $D$ has arclength $2\theta x$, then $%
f(x\,|\,t)=(2\theta x)/(\pi R^2)$; the Law of Cosines gives $\theta $.

\subsection{$R<x<2R$}

\[
f(x)=\dint\limits_{x-R}^R\frac{2x}{\pi R^2}\arccos \left( \frac{t^2+x^2-R^2}{%
2tx}\right) \frac{2t}{R^2}dt 
\]
which corresponds to Parry's (1.18). Straightforward integration provides
the desired result (valid in both of the preceding regions).

\section{Bivariate Case}

We omit geometric details, referring to \cite{Parry} instead. The distance $%
t $ between point $C$ and the origin has density $2t/R^2$ for $0<t<R$. Let $%
f(x,y\,|\,t)$ be the conditional density for distance $x$ between points $B$
and $C$, and distance $y$ between points $A$ and $C$, given $t$. We will
compute the sought-after density $f(x,y)$ for sides $a,b$ via 
\[
f(x,y)=\dint\limits_0^Rf(x,y\,|\,t)\frac{2t}{R^2}dt. 
\]
There are six subcases.

\subsection{$y<x$ and $0<x<R$}

\begin{eqnarray*}
f(x,y) &=&\dint\limits_0^{R-x}\frac{2x}{R^2}\frac{2y}{R^2}\frac{2t}{R^2}%
dt+\dint\limits_{R-x}^{R-y}\frac{2x}{\pi R^2}\arccos \left( \frac{t^2+x^2-R^2%
}{2tx}\right) \frac{2y}{R^2}\frac{2t}{R^2}dt+ \\
&&\dint\limits_{R-y}^R\frac{2x}{\pi R^2}\arccos \left( \frac{t^2+x^2-R^2}{2tx%
}\right) \frac{2y}{\pi R^2}\arccos \left( \frac{t^2+y^2-R^2}{2ty}\right) 
\frac{2t}{R^2}dt
\end{eqnarray*}
which corresponds to formula (4.26) in Parry's thesis \cite{Parry}.

\subsection{$R<x<2R$ and $0<y<2R-x$}

\begin{eqnarray*}
f(x,y) &=&\dint\limits_{x-R}^{R-y}\frac{2x}{\pi R^2}\arccos \left( \frac{%
t^2+x^2-R^2}{2tx}\right) \frac{2y}{R^2}\frac{2t}{R^2}dt+ \\
&&\ \dint\limits_{R-y}^R\frac{2x}{\pi R^2}\arccos \left( \frac{t^2+x^2-R^2}{%
2tx}\right) \frac{2y}{\pi R^2}\arccos \left( \frac{t^2+y^2-R^2}{2ty}\right) 
\frac{2t}{R^2}dt
\end{eqnarray*}
which corresponds to Parry's (4.29). Straightforward integration gives $%
\varphi (x,y)$ (valid in both of the preceding regions).

\subsection{$R<x<2R$ and $2R-x<y<x$}

\[
f(x,y)=\dint\limits_{x-R}^R\frac{2x}{\pi R^2}\arccos \left( \frac{t^2+x^2-R^2%
}{2tx}\right) \frac{2y}{\pi R^2}\arccos \left( \frac{t^2+y^2-R^2}{2ty}%
\right) \frac{2t}{R^2}dt 
\]
which corresponds to Parry's (4.32). This is, of course, $\psi (x,y)$.

\subsection{$x<y$ and $0<y<R$}

\begin{eqnarray*}
f(x,y) &=&\dint\limits_0^{R-y}\frac{2x}{R^2}\frac{2y}{R^2}\frac{2t}{R^2}%
dt+\dint\limits_{R-y}^{R-x}\frac{2x}{R^2}\frac{2y}{\pi R^2}\arccos \left( 
\frac{t^2+y^2-R^2}{2ty}\right) \frac{2t}{R^2}dt+ \\
&&\ \dint\limits_{R-x}^R\frac{2x}{\pi R^2}\arccos \left( \frac{t^2+x^2-R^2}{%
2tx}\right) \frac{2y}{\pi R^2}\arccos \left( \frac{t^2+y^2-R^2}{2ty}\right) 
\frac{2t}{R^2}dt
\end{eqnarray*}
which corresponds to Parry's (4.35). This is, of course, $\varphi (y,x)$.

\subsection{$R<y<2R$ and $0<x<2R-y$}

\begin{eqnarray*}
f(x,y) &=&\dint\limits_{y-R}^{R-x}\frac{2x}{R^2}\frac{2y}{\pi R^2}\arccos
\left( \frac{t^2+y^2-R^2}{2ty}\right) \frac{2t}{R^2}dt+ \\
&&\ \ \dint\limits_{R-x}^R\frac{2x}{\pi R^2}\arccos \left( \frac{t^2+x^2-R^2%
}{2tx}\right) \frac{2y}{\pi R^2}\arccos \left( \frac{t^2+y^2-R^2}{2ty}%
\right) \frac{2t}{R^2}dt
\end{eqnarray*}
which corresponds to Parry's (4.38). This is, of course, $\varphi (y,x)$.

\subsection{$R<y<2R$ and $2R-y<x<y$}

\[
f(x,y)=\dint\limits_{y-R}^R\frac{2x}{\pi R^2}\arccos \left( \frac{t^2+x^2-R^2%
}{2tx}\right) \frac{2y}{\pi R^2}\arccos \left( \frac{t^2+y^2-R^2}{2ty}%
\right) \frac{2t}{R^2}dt
\]
which corresponds to Parry's (4.41). This is, of course, $\psi (y,x)$.

\section{Characteristic Function}

We follow an approach found in \cite{Is1, Is2}. Let $u,v,w$ denote the
squared distances between $A,B,C$ and the origin $O$. Let $\varphi $ denote
the angle between vectors $\overrightarrow{OA},\overrightarrow{OB}$ and $%
\psi $ denote the angle between vectors $\overrightarrow{OA},\overrightarrow{%
OC}$. We have 
\[
a^2=v+w-2\sqrt{vw}\cos (\psi -\varphi ), 
\]
\[
b^2=u+w-2\sqrt{uw}\cos (\psi ), 
\]
\[
c^2=u+v-2\sqrt{uv}\cos (\varphi ) 
\]
by the Law of Cosines, where $u,v,w$ are independent uniform on $[0,R^2]$
and $\varphi ,\psi $ are independent uniform on $[0,2\pi ]$. The
characteristic function for $(a^2,b^2,c^2)$ is thus 
\begin{eqnarray*}
g(r,s,t) &=&\frac 1{R^6}\frac 1{4\pi ^2}\dint\limits_0^{2\pi
}\dint\limits_0^{2\pi }d\varphi \,d\psi
\dint\limits_0^{R^2}\dint\limits_0^{R^2}\dint\limits_0^{R^2}du\,dv\,dw\,\exp
\left[ ir\left( v+w-2\sqrt{vw}\cos (\psi -\varphi )\right) +\right. \\
&&\ \ \ \left. is\left( u+w-2\sqrt{uw}\cos (\psi )\right) +it\left( u+v-2%
\sqrt{uv}\cos (\varphi )\right) \right] .
\end{eqnarray*}
It is well-known that 
\[
\begin{array}{ccc}
\limfunc{E}(c^2)=\dfrac 1i\left. \dfrac{\partial g}{\partial t}\right|
_{r=s=t=0}, &  & \limfunc{E}(b^2\,c^2)=\dfrac 1{i^2}\left. \dfrac{\partial
^2g}{\partial s\,\partial t}\right| _{r=s=t=0}
\end{array}
\]
and the former becomes 
\begin{eqnarray*}
\limfunc{E}(c^2)\ &=&\frac 1i\left. \frac \partial {\partial t}\frac
1{R^4}\frac 1{2\pi }\dint\limits_0^{2\pi }d\varphi
\,\dint\limits_0^{R^2}\dint\limits_0^{R^2}du\,dv\exp \left[ it\left( u+v-2%
\sqrt{uv}\cos (\varphi )\right) \right] \right| _{t=0} \\
\ &=&\frac 1i\left. \frac \partial {\partial t}\frac
1{R^4}\,\dint\limits_0^{R^2}\dint\limits_0^{R^2}\exp \left( it\left(
u+v\right) \right) J_0\left( 2t\sqrt{uv}\right) du\,dv\right| _{t=0}
\end{eqnarray*}
where $J_0(\theta )$ is the zeroth Bessel function of the first kind; hence 
\begin{eqnarray*}
\limfunc{E}(c^2) &=&\frac 1i\left. \frac
1{R^4}\,\dint\limits_0^{R^2}\dint\limits_0^{R^2}\frac \partial {\partial
t}\exp \left( it\left( u+v\right) \right) J_0\left( 2t\sqrt{uv}\right)
\right| _{t=0}du\,dv \\
\ &=&\frac 1i\frac
1{R^4}\,\dint\limits_0^{R^2}\dint\limits_0^{R^2}i(u+v)du\,dv=R^2
\end{eqnarray*}
which is consistent with before. The latter becomes 
\begin{eqnarray*}
\limfunc{E}(b^2\,c^2) &=&\dfrac 1{i^2}\dfrac{\partial ^2}{\partial
s\,\partial t}\frac 1{R^6}\frac 1{4\pi ^2}\dint\limits_0^{2\pi
}\dint\limits_0^{2\pi }d\varphi \,d\psi
\dint\limits_0^{R^2}\dint\limits_0^{R^2}\dint\limits_0^{R^2}du\,dv\,dw\exp
\left[ is\left( u+w-2\sqrt{uw}\cos (\psi )\right) +\right. \\
&&\ \ \ \ \left. \left. it\left( u+v-2\sqrt{uv}\cos (\varphi )\right)
\right] \right| _{s=t=0} \\
\ &=&-\left. \dfrac{\partial ^2}{\partial s\,\partial t}\frac
1{R^6}\dint\limits_0^{R^2}\dint\limits_0^{R^2}\dint\limits_0^{R^2}\exp
\left( is\left( u+w\right) \right) J_0\left( 2s\sqrt{uw}\right) \exp \left(
it\left( u+v\right) \right) J_0\left( 2t\sqrt{uv}\right) du\,dv\,dw\right|
_{s=t=0} \\
\ &=&-\left. \frac
1{R^6}\dint\limits_0^{R^2}\dint\limits_0^{R^2}\dint\limits_0^{R^2}\dfrac{%
\partial ^2}{\partial s\,\partial t}\exp \left( is\left( u+w\right) \right)
J_0\left( 2s\sqrt{uw}\right) \exp \left( it\left( u+v\right) \right)
J_0\left( 2t\sqrt{uv}\right) \right| _{s=t=0}du\,dv\,dw \\
\ &=&-\frac
1{R^6}\dint\limits_0^{R^2}\dint\limits_0^{R^2}\dint%
\limits_0^{R^2}-(u+v)(u+w)du\,dv\,dw=\frac{13}{12}R^4
\end{eqnarray*}
as was to be shown. The fact that $13/12-1=1/12\neq 0$ offers the simplest
proof we know that arbitrary sides of a random triangle in $D$ must be
dependent.

\section{Catalan Numbers}

Let $R=1$ for the remainder of our discussion. From 
\[
\limfunc{P}\left( a^2<x\right) =\limfunc{P}\left( a<\sqrt{x}\right)
=\dint\limits_0^{\sqrt{x}}\left( \dfrac{4\xi }\pi \arccos \left( \dfrac \xi
2\right) -\dfrac{\xi ^2}\pi \sqrt{4-\xi ^2}\right) d\xi 
\]
we obtain that the density for $a^2$ is 
\[
\begin{array}{ccc}
\left( \dfrac{4\sqrt{x}}\pi \arccos \left( \dfrac{\sqrt{x}}2\right) -\dfrac
x\pi \sqrt{4-x}\right) \dfrac 1{2\sqrt{x}}, &  & 0<x<4.
\end{array}
\]
On the one hand, the characteristic function for $a^2$ is

\[
\dint\limits_0^1\dint\limits_0^1\exp \left( it\left( u+v\right) \right)
J_0\left( 2t\sqrt{uv}\right) du\,dv 
\]
by the preceding section; on the other hand, it is 
\begin{eqnarray*}
&&\ \ \dint\limits_0^4\exp (itx)\left( \dfrac{4\sqrt{x}}\pi \arccos \left( 
\dfrac{\sqrt{x}}2\right) -\dfrac x\pi \sqrt{4-x}\right) \dfrac 1{2\sqrt{x}}dx
\\
\ &=&\frac it\left[ 1-\exp (2it)\left( J_0(2t)-iJ_1(2t)\right) \right] \\
\ &=&\frac it\left[ 1-h(t)\right]
\end{eqnarray*}
where $J_1(\theta )=-J_0^{\prime }(\theta )$. A direct evaluation of the
double integral seems to be difficult. Boersma \cite{Boe}, using work of
Zernike \&\ Nijboer \cite{ZN1, ZN2, ZN3}, gave a rapidly-convergent series
for the inner integral: 
\[
\dint\limits_0^1\exp \left( itu\right) J_0\left( 2t\sqrt{uv}\right) du=\frac{%
\sqrt{\pi }}{t^{3/2}v^{1/2}}\exp \left( \frac{it}2\right)
\dsum\limits_{n=0}^\infty (-i)^n(2n+1)J_{n+1/2}\left( \frac t2\right)
J_{2n+1}\left( 2t\sqrt{v}\right) 
\]
but this apparently does not help with the outer integral.

Let $I_0(\theta )$ be the zeroth modified Bessel function of the first kind
and $I_1(\theta )=I_0^{\prime }(\theta )$. We note that the exponential
generating function for the Catalan numbers \cite{Slo}: 
\[
\exp (2t)\left( I_0(2t)-I_1(2t)\right) =\dsum\limits_{n=0}^\infty \frac
1{(n+1)!}\binom{2n}nt^n 
\]
is remarkably similar to the expression for $h(t)$. Replacing $t$ by $it$,
we obtain 
\[
h(t)=\exp (2it)\left( J_0(2t)-iJ_1(2t)\right) =\dsum\limits_{n=0}^\infty
\frac 1{(n+1)!}\binom{2n}n(it)^n 
\]
because $J_0(i\theta )=I_0(\theta )$, $J_1(i\theta )=iI_1(\theta )$.
Therefore the Catalan numbers are associated with the characteristic
function for $a^2$. We wonder if a two-dimensional integer array, suitably
generalizing the Catalan numbers, can be associated with the characteristic
function for $(a^2,b^2)$: 
\[
\dint\limits_0^1\dint\limits_0^1\dint\limits_0^1\exp \left( is\left(
u+w\right) \right) J_0\left( 2s\sqrt{uw}\right) \exp \left( it\left(
u+v\right) \right) J_0\left( 2t\sqrt{uv}\right) du\,dv\,dw. 
\]
Since the bivariate density $f(x,y)$ for $(a,b)$ is much more complicated
than the univariate density $f(x)$ for $a$, an answer to our question may be
a long time coming.

\section{Acknowledgement}

I am thankful to Michelle Parry for her correspondence. Much more relevant
material can be found at \cite{F1, F2}, including experimental computer runs
that aided theoretical discussion here.


\begin{thebibliography}{99}
\bibitem{Delt}  R. Deltheil, \textit{Probabilit\'es G\'eom\'etriques}, t. 2, 
\textit{Trait\'e du calcul des Probabili}$J$\textit{t\'es et de ses
Applications}, f. 2, ed E. Borel, Gauthier-Villars, 1926, pp. 40--42,
114--120.

\bibitem{Ham}  J. M. Hammersley, The distribution of distance in a
hypersphere, \textit{Annals Math. Statist.} 21 (1950) 447--452; MR0037481
(12,268e).

\bibitem{Lord}  R. D. Lord, The distribution of distance in a hypersphere, 
\textit{Annals Math. Statist.} 25 (1954) 794--798; MR0065048 (16,377d).

\bibitem{Alag}  V. S. Alagar, The distribution of the distance between
random points, \textit{J. Appl. Probab.} 13 (1976) 558--566; MR0418183 (54
2\#6225).

\bibitem{Solo}  H. Solomon, \textit{Geometric Probability}, SIAM, 1978, pp.
35--36, 128--129; MR0488215 (58 \#7777).

\bibitem{Dunb}  S. R. Dunbar, The average distance between points in
geometric figures, \textit{College Math. J.} 28 (1997) 187--197; MR1444006
(98a:52007).

\bibitem{TuFi}  S.-J. Tu and E. Fischbach, Random distance distribution for
spherical objects: general theory and applications to physics, \textit{J.
Phys. A} 35 (2002) 6557--6570; MR1928848.

\bibitem{Parry}  M. Parry, \textit{Application of Geometric Probability
Techniques to Elementary Particle and Nuclear Physics}, Ph.D. thesis, Purdue
Univ., 1998.

\bibitem{PaFi}  M. Parry and E. Fischbach, Probability distribution of
distance in a uniform ellipsoid: theory and applications to physics, \textit{%
J. Math. Phys.} 41 (2000) 2417--2433; MR1751899 (2001j:81267).

\bibitem{Is1}  Y. Isokawa, Limit distributions of random triangles in
hyperbolic planes, \textit{Bull. Faculty Educ. Kagoshima Univ. Natur. Sci.}
49 (1997) 1--16; MR1653095 (99k:60020).

\bibitem{Is2}  Y. Isokawa, Geometric probabilities concerning large random
triangles in the hyperbolic plane, \textit{Kodai Math. J.} 23 (2000)
171--186; MR1768179 (2001f:60014).

\bibitem{Boe}  J. Boersma, On the computation of Lommel's functions of two
variables, \textit{Math. Comp.} 16 (1962) 232--238; MR0146419 (26 \#3941).

\bibitem{ZN1}  F. Zernike and B. R. A. Nijboer, Th\'eorie de la diffraction
des aberrations, \textit{La Th\'eorie des Images Optiques}, Proc. 1946 Paris
colloq., ed. P. Fleury, A. Mar\'echal and C. Anglade, La Revue d'Optique,
1949, pp. 227--235.

\bibitem{ZN2}  B. R. A. Nijboer, \textit{The Diffraction Theory of
Aberrations}, Ph.D. thesis, Univ. of Groningen, 1942, available online at
http://www.nijboerzernike.nl/\_html/intro.html.

\bibitem{ZN3}  A. J. E. M. Janssen, J. J. M. Braat and P. Dirksen, On the
computation of the Nijboer-Zernike aberration integrals at arbitrary
defocus, \textit{J. Mod. Optics} 51 (2004) 687--703; available online at
http://www.nijboerzernike.nl/\_html/biblio.html.

\bibitem{Slo}  N. J. A. Sloane, On-Line Encyclopedia of Integer Sequences,
A144186.

\bibitem{F1}  S. Finch, Random triangles. I--IV, unpublished essays (2010),
http://algo.inria.fr/bsolve/.

\bibitem{F2}  S. Finch, Simulations in R\ involving triangles and
tetrahedra, http://algo.inria.fr/csolve/rsimul.html.\\

\begin{tabular}{lll}
& Steven Finch &  \\ 
& Dept. of Statistics &  \\ 
& Harvard University &  \\ 
& Cambridge, MA, USA &  \\ 
& \textit{Steven.Finch@inria.fr} & 
\end{tabular}
\end{thebibliography}
\end{document}